\documentclass[12pt]{article}

\usepackage{amsfonts}
\usepackage{latexsym}
\usepackage{eepic}

\def\baselinestretch{1} \topmargin -12pt \headsep 0pt \footskip 30pt
\def\medskipamount{12pt}

\begin{document}


\newcommand{\bit}[1]{\pagebreak[3] \section{#1} \setcounter{equation}{0}}
\renewcommand{\theequation}{\thesection .\arabic{equation}}


\newcommand{\re}[1]{\mbox{\bf (\ref{#1})}}


\catcode`\@=\active
\catcode`\@=11
\def\@eqnnum{\hbox to .01pt{}\rlap{\bf \hskip -\displaywidth(\theequation)}}
\catcode`\@=12



\catcode`\@=\active
\catcode`\@=11
\newcommand{\nc}{\newcommand}

   
\nc{\bs}[1]{\addvspace{\medskipamount} \pagebreak[3]
\refstepcounter{equation}
\noindent \begin{em}{\bf (\theequation) #1.} \nopagebreak}

\nc{\es}{\par \end{em} \addvspace{\medskipamount} } 

\nc{\br}[1]{\addvspace{\medskipamount} \pagebreak[3]
\refstepcounter{equation} 
\noindent {\bf (\theequation) #1.} \nopagebreak }

\nc{\er}{\par \addvspace{\medskipamount} }

\newcounter{index}
\nc{\bl}{\begin{list}{{\rm (\roman{index})}}{\usecounter{index}}}
   
\nc{\el}{\end{list}}

\nc{\pf}{\addvspace{\medskipamount} \par \noindent {\em Proof}}

\nc{\fp}{\phantom{.} \hfill \mbox{     $\Box$} \par
\addvspace{\medskipamount}} 


\nc{\C}{\mathbb C}
\nc{\N}{\mathbb N}
\nc{\Pj}{\mathbb P}
\nc{\Q}{\mathbb Q}
\nc{\R}{\mathbb R}
\nc{\Z}{\mathbb Z}
 

\nc{\oper}[1]{\mathop{\mathchoice{\mbox{\rm #1}}{\mbox{\rm #1}}
{\mbox{\rm \scriptsize #1}}{\mbox{\rm \tiny #1}}}\nolimits}
\nc{\Aut}{\oper{Aut}}
\nc{\ev}{\oper{ev}}
\nc{\Gr}{\oper{Gr}}
\nc{\Hom}{\oper{Hom}}
\nc{\itHom}{\oper{{\it Hom}}}
\nc{\Jac}{\oper{Jac}}
\nc{\seg}{\oper{{\it s}}}
\nc{\Sym}{\oper{Sym}}

\nc{\red}{{\rm red}}
\nc{\vir}{{\rm vir}}

\nc{\operlimits}[1]{\mathop{\mathchoice{\mbox{\rm #1}}{\mbox{\rm #1}}
{\mbox{\rm \scriptsize #1}}{\mbox{\rm \tiny #1}}}}

\nc{\Coeff}{\operlimits{Coeff}}
\nc{\prodstar}{\operlimits{${\displaystyle \prod^{}_{}}^*$}}


\nc{\al}{\alpha}
\nc{\ga}{\gamma}
\nc{\la}{\lambda}
\nc{\La}{\Lambda}
\nc{\si}{\sigma}


\nc{\Si}{C}
\nc{\curve}{\Sigma}


\nc{\bigdownarg}[1]{{\phantom{\scriptstyle #1}\Big\downarrow
\raisebox{.4ex}{$\scriptstyle #1$}}} 
\nc{\leftbigdownarg}[1]{{\raisebox{.4ex}{$\scriptstyle #1$}
\Big\downarrow \phantom{\scriptstyle #1}}} 
\nc{\bigdowneq}{{\Big| \! \Big|}}

\nc{\cl}{{\mathcal L}}
\nc{\co}{{\mathcal O}}

\nc{\GL}[1]{{{\rm GL(}#1{\rm )}}}
\nc{\Sp}[1]{{{\rm Sp(}#1{\rm )}}}

\nc{\bino}[2]{\mbox{\Large $#1 \choose #2$}}

\nc{\half}{\mathchoice{{\textstyle \frac{\scriptstyle 1}{\scriptstyle
2}}} {{\textstyle\frac{\scriptstyle 1}{\scriptstyle 2}}}
{\frac{\scriptscriptstyle 1}{\scriptscriptstyle 2}}
{\frac{\scriptscriptstyle 1}{\scriptscriptstyle 2}}}
\nc{\quarter}{\mathchoice{{\textstyle \frac{\scriptstyle
1}{\scriptstyle 4}}} {{\textstyle\frac{\scriptstyle 1}{\scriptstyle
4}}} {\frac{\scriptscriptstyle 1}{\scriptscriptstyle 4}}
{\frac{\scriptscriptstyle 1}{\scriptscriptstyle 4}}}
\nc{\ratio}[2]{\mathchoice{ {\textstyle \frac{\scriptstyle
#1}{\scriptstyle #2}}} {{\textstyle\frac{\scriptstyle #1}{\scriptstyle
#2}}} {\frac{\scriptscriptstyle #1}{\scriptscriptstyle #2}}
{\frac{\scriptscriptstyle #1}{\scriptscriptstyle #2}}}

\nc{\lrow}{\longrightarrow}
\nc{\mbar}{\overline{M}}
\nc{\sans}{\backslash}
\nc{\st}{\, | \,}

\catcode`\@=12



\noindent
{\LARGE \bf On the quantum cohomology}\smallskip \\
{\LARGE \bf of a symmetric product 
of an algebraic curve}
\bigskip \\ 
{\bf Aaron Bertram } \smallskip \\
Department of Mathematics, University of Utah,\\
Salt Lake City, Utah 84112 USA \medskip \\
{\bf Michael Thaddeus } \smallskip \\
Department of Mathematics, Columbia University, \\
New York, N.Y. 10027 USA 
\renewcommand{\thefootnote}{}
\footnotetext{A.B. supported by NSF grant DMS--95--00865;
M.T. supported by NSF grant DMS--95--00964.}

\medskip
\smallskip

Quantum cohomology is a novel multiplication on the cohomology of a
smooth complex projective variety, or even a compact symplectic
manifold.  It can be regarded as a deformation of the ordinary cup
product, defined in terms of the Gromov-Witten invariants of the
manifold.  Since its introduction in 1991, there has been enormous
interest in computing quantum cohomology for various target spaces.
Attention has focused on homogeneous spaces, complete intersections,
surfaces, and of course on Calabi-Yau threefolds, where it is a key
part of mirror symmetry.  However, to the authors' knowledge, no one
has yet studied the quantum cohomology of a symmetric product of a
smooth curve.  This is strange, because the problem is attractive from
several points of view.

First, though quantum cohomology is clearly a fundamental invariant of
a variety or a symplectic manifold, it has been explicitly
computed in relatively few cases: just some of those mentioned above.

A more specific motivation comes from a link with Seiberg-Witten
theory.  For a complex curve $\Si$, the Seiberg-Witten-Floer
cohomology of the real 3-manifold $\Si \times S^1$ is isomorphic to
the cohomology of the $d$th symmetric product $\Si_d$
\cite{d,mst,mun}, and is expected to carry a natural product that
corresponds to the quantum product \cite{cwzm,mar,pss}.  The index $d$
depends on the spin-c structure on $\Si \times S^1$ chosen in
Seiberg-Witten theory.  In the ordinary Seiberg-Witten theory, this
isomorphism only holds when $d<g-1$, where $g$ is the genus of $C$.
However, higher symmetric products enter the picture if the
Seiberg-Witten functional is perturbed in the manner of Taubes
\cite{t} and Morgan-Szab\'o-Taubes \cite{mst}.

The quantum cohomology of a symmetric product is therefore the base
ring of the ``quantum category'' in Seiberg-Witten theory, as
introduced by Donaldson \cite[\S10.4]{ms}.

On the other hand, there is also a link with Brill-Noether theory and
the stratification of $\Si_d$ by the special linear series.  This
makes it feasible to compute the quantum cohomology using algebraic
geometry.  To carry out the necessary enumerative computations on the
strata, crucial use is made of a generalization of Porteous's formula,
due to Harris and Tu \cite{ht}.

There are several kinds of quantum cohomology; this paper is concerned
only with the so-called ``quantum cohomology algebra'' \cite{g} or
``little quantum cohomology'' \cite{fp} of $\Si_d$.  This is a ring
generated by the rational cohomology together with one deformation
parameter for each generator of $H_2(\Si_d; \Z)$.  Its structure
constants are given by 3-point Gromov-Witten invariants.  Not quite
all of these invariants are calculated herein, but in a sense most of
them are.  More precisely, the following results are proved.

\begin{itemize}
  
\item The number of possible deformation parameters is the second
  Betti number of $\Si_d$, which is fairly large.  Nevertheless, the
  quantum product is shown to depend nontrivially only on a
  single parameter $q$ (Proposition \ref{z}).
  
\item Explicit formulas are given for the coefficients of $q$
  (Corollary \ref{t}) and $q^2$ (Corollary \ref{dd}) in the quantum
  product.

\item All the terms in the quantum product are computed for
  $\Si_{g-1}$ (Corollary \ref{aa} ii).
  
\item The coefficient of $q^e$ in the quantum product for $\Si_d$ is
  shown to vanish if $d<g-1$ and $e > \frac{d-3}{g-1-d}$ (Corollary
  \ref{bb}), or if $d>g-1$ and $e>1$ (Corollary \ref{aa} i).

\end{itemize}

Schematically, the coefficients vanish in the regions marked A and B
in the diagram below.  A follows from a straightforward dimension
count; B is more subtle.  Only in the regions under the hyperbola, or
under the line $e=1$, are there nonzero coefficients.

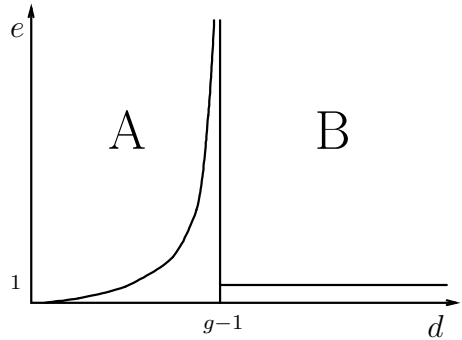
\setlength{\unitlength}{0.00041666in}
\begingroup\makeatletter\ifx\SetFigFont\undefined
\def\x#1#2#3#4#5#6#7\relax{\def\x{#1#2#3#4#5#6}}%
\expandafter\x\fmtname xxxxxx\relax \def\y{splain}%
\ifx\x\y   
\gdef\SetFigFont#1#2#3{%
  \ifnum #1<17\tiny\else \ifnum #1<20\small\else
  \ifnum #1<24\normalsize\else \ifnum #1<29\large\else
  \ifnum #1<34\Large\else \ifnum #1<41\LARGE\else
     \huge\fi\fi\fi\fi\fi\fi
  \csname #3\endcsname}%
\else
\gdef\SetFigFont#1#2#3{\begingroup
  \count@#1\relax \ifnum 25<\count@\count@25\fi
  \def\x{\endgroup\@setsize\SetFigFont{#2pt}}%
  \expandafter\x
    \csname \romannumeral\the\count@ pt\expandafter\endcsname
    \csname @\romannumeral\the\count@ pt\endcsname
  \csname #3\endcsname}%
\fi
\fi\endgroup
\begin{center}
\begin{picture}(5742,4113)(0,-10)
\thicklines
\path(330.000,3966.000)(300.000,4086.000)(270.000,3966.000)
\path(300,4086)(300,336)(5700,336)
\path(5580.000,306.000)(5700.000,336.000)(5580.000,366.000)
\path(2700,336)(2700,3926)
\path(2700,561)(5580,561)

\path(475.856,336)
        (526.219,340.734)
        (572.820,346.042)
        (615.879,351.055)
        (655.616,355.801)
        (692.250,360.306)
        (757.090,368.707)
        (812.156,376.476)
        (859.207,383.834)
        (900.000,391.000)

\path(900,391)  (960.242,402.580)
        (1033.658,417.424)
        (1073.906,425.880)
        (1115.769,434.924)
        (1158.685,444.480)
        (1202.096,454.473)
        (1245.442,464.825)
        (1288.162,475.462)
        (1329.696,486.307)
        (1369.486,497.284)
        (1441.592,519.333)
        (1500.000,541.000)

\path(1500,541) (1564.785,570.414)
        (1601.959,588.565)
        (1641.562,608.689)
        (1683.030,630.545)
        (1725.795,653.890)
        (1769.288,678.483)
        (1812.945,704.080)
        (1856.197,730.441)
        (1898.477,757.322)
        (1939.219,784.483)
        (1977.855,811.681)
        (2046.541,865.219)
        (2100.000,916.000)

\path(2100,916) (2143.193,970.966)
        (2187.720,1040.474)
        (2210.020,1079.303)
        (2232.099,1120.114)
        (2253.771,1162.356)
        (2274.851,1205.478)
        (2295.155,1248.928)
        (2314.497,1292.157)
        (2332.693,1334.612)
        (2349.557,1375.743)
        (2364.905,1414.998)
        (2378.551,1451.826)
        (2400.000,1516.000)

\path(2400,1516)        (2415.014,1571.419)
        (2429.106,1634.628)
        (2442.311,1704.517)
        (2448.593,1741.621)
        (2454.668,1779.977)
        (2460.539,1819.450)
        (2466.211,1859.898)
        (2471.689,1901.184)
        (2476.978,1943.170)
        (2482.082,1985.715)
        (2487.006,2028.682)
        (2491.753,2071.931)
        (2496.330,2115.325)
        (2500.740,2158.724)
        (2504.988,2201.989)
        (2509.078,2244.983)
        (2513.016,2287.565)
        (2516.805,2329.597)
        (2520.450,2370.942)
        (2523.956,2411.459)
        (2527.327,2451.010)
        (2530.569,2489.457)
        (2533.685,2526.660)
        (2539.558,2596.781)
        (2544.984,2660.265)
        (2550.000,2716.000)

\path(2550,2716)        (2553.596,2755.584)
        (2557.209,2797.839)
        (2560.867,2843.204)
        (2564.597,2892.121)
        (2568.427,2945.027)
        (2572.385,3002.362)
        (2576.497,3064.566)
        (2580.791,3132.078)
        (2585.295,3205.337)
        (2587.634,3244.259)
        (2590.036,3284.784)
        (2592.504,3326.964)
        (2595.042,3370.857)
        (2597.652,3416.515)
        (2600.340,3463.996)
        (2603.106,3513.352)
        (2605.957,3564.640)
        (2608.894,3617.914)
        (2611.921,3673.229)
        (2615.042,3730.640)
        (2618.259,3790.203)
        (2621.578,3851.971)
        (2625.000,3926.000)

\put(0,3711){\makebox(0,0)[lb]{\smash{{{\SetFigFont{12}{14.4}{rm}$e$}}}}}
\put(5300,-100){\makebox(0,0)[lb]{\smash{{{\SetFigFont{12}{14.4}{rm}$d$}}}}}
\put(0,500){\makebox(0,0)[lb]{\smash{{{\SetFigFont{12}{14.4}{rm}$\scriptstyle 1$}}}}}
\put(2450,0){\makebox(0,0)[lb]{\smash{{{\SetFigFont{12}{14.4}{rm}$\scriptstyle g-1$}}}}}
\put(1275,2286){\makebox(0,0)[lb]{\smash{{{\SetFigFont{20}{14.4}{rm}A}}}}}
\put(3900,2286){\makebox(0,0)[lb]{\smash{{{\SetFigFont{20}{14.4}{rm}B}}}}}
\end{picture}
\end{center}

Putting these results together completely determines the quantum
product on $\Si_d$ in all cases except $d \in [\ratio{3}{4}g, g-1)$.

Instead of seeking to characterize the quantum product completely,
one can alternatively ask just for generators and relations for the
quantum ring.  This is a weaker question, because the quantum
relations do not determine the additive isomorphism between the
quantum and classical rings.  It is answered in all cases except $d
\in [\ratio{4}{5}g-\ratio{3}{5}, g-1)$, even using only the first order
terms (Proposition \ref{w} and the discussion following).  In some
cases, one encounters the curious fact that the quantum and classical
rings are isomorphic, but only by an automorphism which is not the
identity.

The organization of the paper is straightforward.  Section 1 recalls
some basic facts on quantum cohomology.  Section 2 recalls some basic
facts on symmetric products of a curve, shows why only one deformation
parameter is involved, and explains the vanishing in region A.\@
Section 3 introduces the Brill-Noether methods which will be used to
compute the 3-point Gromov-Witten invariants.  In particular, the
result of Harris and Tu mentioned above is reviewed.  Using this, the
degree 1 invariants are computed in section 4, and the degree 2
invariants in section 5.  Section 6 is devoted to the vanishing of
higher-degree invariants in region B, and to the computation for
$d=g-1$.  Section 7 explains how to find generators and relations for
the quantum ring.  Finally, sections 8 and 9 are essentially
appendices: section 8 outlines a remarkable connection with Givental's
work on the rational curves on a quintic threefold, while section 9 is
concerned with the first two homotopy groups of $\Si_d$, explaining
Proposition \ref{z} from the point of view of symplectic topology.

A few conventions: cohomology of a space is with rational coefficients
unless otherwise mentioned, and cohomology of a sheaf is over the
curve $\Si$ unless otherwise mentioned.  If $n!$ appears in the
denominator of some expression and $n < 0$, this means that the whole
expression vanishes: that is, $1/n! = 0$ for $n<0$.  

We wish to thank Arnaud Beauville, Olivier Debarre, Tom Graber,
Michael Hutchings, Rahul Pandharipande, Michael Roth, Bernd Siebert,
Ravi Vakil, and Angelo Vistoli for very helpful conversations and
advice.  We also wish to thank the Institut Mittag-Leffler for its
warm hospitality and extraordinary atmosphere, which inspired the
present work.

\medskip

\bit{Preliminaries on quantum cohomology}

Let $X$ be a smooth complex projective variety with fixed
polarization.  Witten \cite{w} introduced two rings associated to $X$,
the {\em big} and {\em little} quantum cohomology.  The big quantum
cohomology was used by Kontsevich \cite{k} to count rational curves in
the plane.  However, there are very few spaces for which it has been
characterized in full.  Most work, including that of Givental and Kim
\cite{g,gk}, and the whole of the present paper, is concerned with the
little quantum cohomology.

\br{Definition}
A {\em stable map} is a morphism $\phi$ from a complete nodal curve
$\curve$ to $X$ such that $\{ \psi \in \Aut \curve \st \phi \psi =
\phi \}$ is finite.  An {\em $n$-pointed stable map} is similar,
except that $n$ distinct smooth points $x_1, \dots, x_n \in \curve$
are also chosen, and the automorphisms $\psi$ are required to fix each
$x_i$.  
\er

The arithmetic genus $g= p_a(\curve)$ and the homology class $e =
\phi_*[\curve] \in H_2(X;\Z)$ are discrete invariants of a stable map.
When these are fixed, a fundamental theorem asserts the existence of
projective moduli spaces $\mbar_g(X,e)$ of stable maps and
$\mbar_{g,n}(X,e)$ of $n$-pointed stable maps, as well as a forgetful
morphism $f_e: \mbar_{g,n}(X,e) \to \mbar_g(X,e)$ and an evaluation
morphism $\ev_e: \mbar_{g,n}(X,e) \to X^n$.  See Fulton-Pandharipande
\cite{fp} for details.  The subscript $e$ will be suppressed when
there is no danger of confusion.  Also, this paper will be concerned
only with the case $g=0$, $n=1$.

It is most accurate to regard these moduli spaces as stacks rather
than schemes.  They are stratified by smooth substacks on which the
dimension of the deformation space is constant.  Since each stable map
has finitely many automorphisms, these strata are Deligne-Mumford
stacks.

This paper, however, adopts a more na\"\i ve point of view, regarding
the moduli spaces as schemes stratified by subschemes with a natural
orbifold structure.  As pointed out by the referee, this is a bit
tricky, since smooth Deligne-Mumford stacks need not be orbifolds in
the optimal algebraic sense of being locally a quotient of an affine
scheme by a finite group.  So strictly speaking, the word {\em
orbifold\/} should henceforth be taken to mean a smooth
Deligne-Mumford stack, and an {\em orbifold vector bundle} should mean
a vector bundle over such a stack.  Fortunately, in the only place
where an explicit calculation with orbifolds is carried out, namely
the proof of \re{n}, the space in question is an orbifold in the
optimal sense mentioned above.

The moduli space $\mbar_0(X,e)$ has expected dimension $\dim X - 3 +
c_1(TX) \cdot e$.  It is not generally of this dimension, or even
equidimensional.  But it is endowed with a natural {\em virtual
fundamental class}, an equivalence class of algebraic cycles
$[\mbar_0(X,e)]^\vir \in A_*(\mbar_0(X,e))$ having the expected
dimension.  The general construction of this virtual class involves
the deformation theory of stable maps and is rather complicated.  It
has been carried out by several authors \cite{bf,lt,s}.  We will use
only the following three basic facts.

\bs{Proposition}
\label{h}
\bl
\item If $T^1_X(\phi)$, $T^2_X(\phi)$ are the first-order deformation
and obstruction spaces of the map $\phi$, and $T^0(\curve)$,
$T^1(\curve)$ are the first-order endomorphism and deformation
spaces of the curve $\curve$, then there is a natural exact sequence
$$
\begin{array}{crclcc}
 & 0 & \lrow & T^0(\curve) & \lrow & H^0(\curve, \phi^* TX) \\
 \lrow & T^1_X(\phi) & \lrow & T^1(\curve) & \lrow 
& H^1(\curve, \phi^* TX) \\ 
 \lrow & T^2_X(\phi) & \lrow & 0. & & 
\end{array}
$$

\item On any reduced locus where the first-order obstruction spaces of
the map have constant dimension, the virtual class is the Euler
class of the orbifold vector bundle formed by these spaces.  In
particular, if the first-order obstructions vanish, then the virtual
class is simply the orbifold fundamental class.

\item The forgetful morphism $f$ is flat, and 
$$[\mbar_{0,n}(X,e)]^\vir = f^*[\mbar_0(X,e)]^\vir.$$   \el \es

\pf. See the work of Ran \cite{ran},  Li-Tian \cite{lt}, 
Behrend \cite{b}, Behrend-Fantechi \cite{bf} and Behrend-Manin
\cite{bm}.  \fp

\br{Example}
\label{j}
If $X = \Pj^r$, then from the long exact sequence of
$$ 0 \lrow \co \lrow \co(1)^{r+1} \lrow T \Pj^r \lrow 0,$$ 
it follows that $H^1(\curve, \phi^* T \Pj^r) = 0$ for a curve $\curve$
of arithmetic genus $0$. Hence the obstructions vanish and
$\mbar_0(\Pj^r, e)$ is an orbifold of the expected dimension $r - 3 +
e(r+1)$.  
\er

\br{Definition}
For $a_1, \dots, a_n \in H^*(X;\Q)$, the
{\em $n$-point degree $e$ Gromov-Witten invariant} is defined as
$$\langle a_1, \dots, a_n \rangle_e = 
\ev_e^* 
{\Big (} \prod_i \pi_i^*
a_i {\Big )}
[\mbar_{0,n}(X,e)]^\vir,$$ 
where $\pi_i: X^n \to X$ is projection on the $i$th factor.
\er

In practice, we will always evaluate Gromov-Witten invariants using
the following equivalent formulation.  Suppose that $\tilde{M}_\al$
are a collection of orbifold resolutions of the closures of the smooth
strata $M_\al$ of $\mbar_0(X,e)$ such that $\tilde{M}_\al
\times_{\mbar_0(X,e)} \mbar_{0,1}(X,e)$ are also orbifold resolutions
of the strata of $\mbar_{0,1}(X,e)$.  Let $\tilde{f}_\al$ and
$\widetilde{\ev}_\al$ be the liftings of $f_e$ and $\ev_e$ to these
fibered products.  Choose a cycle representing the virtual class of
$\mbar_0(X,e)$, and let $[\tilde{M}_\al]^\vir$ be the cycle whose
projection to the closure $\overline{M}_\al$ consists of all
components of the virtual cycle supported in $\overline{M}_\al$ but
not in the closure of any smaller stratum.  The rational equivalence
classes of the individual cycles $[\tilde{M}_\al]^\vir$ might depend
on the choice of the representative, but this will not impair our
arguments.

\bs{Proposition}
\label{ll}
With the above notation,
$$\langle a_1, \dots, a_n \rangle_e = 
\sum_\al {\Big (} \prod_i (\tilde{f}_\al)_*
\, \widetilde{\ev}_\al^* \, a_i {\Big )} [\tilde{M}_\al]^\vir.$$
In particular, if $\mbar_0(X,e)$ is an orbifold of the expected
dimension,
$$\langle a_1, \dots, a_n \rangle_e = 
{\Big (} \prod_i f_*
\, \ev^* a_i {\Big )} [\mbar_0(X,e)].$$
\es

\pf.  First of all, the forgetful morphism $f_e$ factors through a
birational morphism to the $n$th fibered power
${\mbar_{0,1}(X,e)}^n_{\mbar_0(X,e)}$, which is also flat over
$\mbar_0(X,e)$.  An evaluation map to $X^n$ is still defined on this
space, so it suffices to perform the computation here.

It is convenient to denote the disjoint
union of the $\tilde{M}_\al$ by $\tilde{M}_0$, its fibered
product over $\mbar_0(X,e)$ with $\mbar_{0,1}(X,e)$ by
$\tilde{M}_{0,1}$, and the forgetful and evaluation maps on
$\tilde{M}_{0,1}$ by $\tilde{f}$ and $\widetilde{\ev}$.  Then the
proposition follows immediately from \re{h}(iii) and the fact that the
diagram
$$
\renewcommand{\arraystretch}{1.5}
\begin{array}{ccccc}
{(\tilde{M}_{0,1})}^n_{\tilde{M}_0} & \lrow &
\tilde{M}_{0,1}^n & \stackrel{\widetilde{\ev}}{\lrow} & X^n \\
\leftbigdownarg{\tilde{f}} & & \bigdownarg{\tilde{f}} & & \\
\tilde{M}_0 & \lrow &
\tilde{M}_0^n & & 
\end{array}$$
is a fiber square, where the top and bottom of the square
are diagonal embeddings.  \fp

A fundamental result states that Gromov-Witten invariants
are deformation invariants \cite{bf,fo,lt}.
Also, they satisfy the following properties. 

\bs{Lemma}
\label{d}
\bl
\item Gromov-Witten invariants are symmetric on classes of even
  degree, antisymmetric on classes of odd degree.

\item If $e \neq 0$ and $a_1 \in H^0$ or $H^1$, then $\langle a_1,
\dots, a_n \rangle_e = 0$.

\item If $e \neq 0$ and $a_1 \in H^2$, then 
$\langle a_1, \dots, a_n \rangle_e = (a_1 \cdot e) \langle a_2,
\dots, a_n \rangle_e$.

\el
\es

\pf.  See for example Behrend \cite{b} and Behrend-Manin \cite{bm}.  \fp

The quantum product has coefficients in the following ring.

\br{Definition}
The {\em Novikov ring} $\La$ is the subring of $\Q[[H_2(X;\Z)]]$
consisting of formal power series $\sum_{e \in H_2} \la_e q^e$, where
$q$ is a formal variable, $\la_e \in \Q$, and 
$$\{ e \in H_2 \st \la_e \neq 0, \omega \cdot e \leq \kappa \}$$
is finite for all $\kappa \in \Q$, $\omega$ being the class of a fixed
polarization on $X$.  
\er

As an example, if $H_2 \cong \Z$, then $\La \cong \Q[[q]][q^{-1}]$.

One extends the Gromov-Witten invariants linearly to $H^*(X;\La)$.

\br{Definition}
Let $\ga_i$ be a basis for $H^*(X;\Q)$ and $\ga^i$ the dual basis
with respect to the Poincar\'e pairing.  For $a,b \in H^*(X;\La)$, the
{\em little quantum product} is
$$ a * b = \sum_e \sum_i q^e \ga_i \langle a,b,\ga^i \rangle_e.$$
\er

Alternatively, one can define the product without choosing a basis or
mentioning the Gromov-Witten invariants explicitly, as
$$a * b = \sum_e q^e \ev_* f^* \bigl( (f_* \ev^* a) \cdot (f_* \ev^* b)
\bigr),$$
where $f$ is restricted to the virtual cycle.  This reveals
that the cup product has in some sense been transferred from $X$ to
the moduli space of stable maps.

Since $\mbar_{0,1}(X,0) = \mbar_0(X,0) = X$, the quantum product
equals the cup product modulo $q$.  Furthermore, it follows from 
\re{d}(ii) that the quantum product with any element of $H^1$ equals
the cup product.

A fundamental theorem asserts that the little quantum product is
associative \cite{bf,fo,lt,rt}.

\br{Definition}
The {\em little quantum cohomology} $QH^*(X)$ is defined to be the ring
additively isomorphic to $H^*(X;\La)$, but with the little quantum
product as multiplication.  
\er

The expected complex dimension of $\mbar_0(X,e)$ is $\dim X - 3 +
c_1(TX) \cdot e$, so $QH^*(X)$ is graded if $q^e$ is given degree $2
c_1(TX) \cdot e$.

\bit{Symmetric products of a curve}

We begin by recalling some basic facts on the cohomology of symmetric
products of a curve.  Good references are Arbarello et al.\
\cite{acgh} and Macdonald \cite{mac}.

Let $\Si$ be a smooth projective curve of genus $g$, and let $\Si_d$
be the $d$th symmetric product, which is a smooth projective variety
of dimension $d$.  It can be regarded as the moduli space of effective
divisors $D$ on $\Si$ of degree $|D| = d$.  Accordingly, there exists
a universal divisor $\Delta \subset \Si_d \times \Si$ having the
obvious property.

The Poincar\'e dual of $\Delta$ is a class in $H^2(\Si_d \times \Si;
\Z)$, which determines a map $\mu: H_*(\Si; \Z) \to H^*(\Si_d; \Z)$.
Let $e_i$ be a basis of $H_1(\Si; \Z)$ for which the
intersection form is the standard symplectic form. 
Then define $\xi_i = \mu(e_i) \in H^1(\Si_d)$ and $\eta = \mu(1) \in
H^2(\Si_d)$.  It is also convenient to define $\si_i = \xi_i
\xi_{i+g}$ for $i \leq g$ and $\theta = \sum_{i=1}^g \si_i$.

There is a morphism $\iota: \Si_{d-1} \to \Si_d$ given by $D
\mapsto D + p$; it is an embedding, and its image is a divisor in
$\Si_d$.  Indeed, this divisor is exactly $\Delta|_{\Si_d \times p}$,
so it is Poincar\'e dual to $\eta$.

There is also a natural morphism $AJ: \Si_d \to \Jac_d(\Si)$, the {\em
Abel-Jacobi map}, taking a divisor $D$ to the line bundle $\co(D)$.
Since every nonzero section of a line bundle determines a divisor, and
vice versa up to scalars, every fiber $AJ^{-1}(L)$ is a projective
space, namely $\Pj H^0(L)$.

Indeed, for $d \geq 2g-2$, $\Si_d$ is the projectivization of a vector
bundle over $\Jac_d$, as follows.  Fix $p \in \Si$, and let $\cl$ be a
Poincar\'e line bundle over $\Jac_d \times \Si$, normalized so that
$\cl|_{\Jac_d \times p} = \co$.  Let $\pi: \Jac_d \times \Si \to
\Jac_d$ be the projection, and let $U = \pi_* \cl$.  Then $U$ is a
vector bundle of rank $d-g+1$ by Riemann-Roch, and there is a
natural isomorphism $\Si_d = \Pj U$. 

If $\co(1)$ is the twisting sheaf of this projective bundle, the
embedding $\iota$ satisfies $\iota^*\co(1) = \co(1)$ for any $d >
2g-2$. So for $d < 2g-2$, it is reasonable to define $\co(1) = \iota^*
\co(1)$ by descending induction.  For any $d \geq 0$ and any $D \in
\Si_d$, the fiber of $\co(-1)$ at $D$ is then naturally isomorphic to
the space of sections of $\cl|_{AJ(D) \times \Si}$ vanishing at $D$,
where $\cl$ is again a Poincar\'e bundle normalized as above.  To
construct the isomorphism, one simply multiplies by the appropriate
power of the natural section of $\co(p)$ vanishing at $p$.

\bs{Proposition} 
\label{x}
\bl
\item The twisting sheaf $\co(1)$ of $\Pj U$ is isomorphic to
$\co(\Si_{d-1})$, 
so $c_1(\co(1)) = \eta$. 

\item There is a natural isomorphism $\cl(1) = \co(\Delta)$, where
  $\cl$ is short for $AJ^* \cl$.

\item If $e_i$ is viewed as an element of $H^1(\Jac_d) \cong H_1(\Si)$,
then $\xi_i = AJ^* e_i$.

\item If $\Theta$ is a theta-divisor on $\Jac_d$, then $AJ^*
c_1(\co(\Theta)) = \theta$.
\el
\es

\pf. On $\Pj U = \Si_d$, evaluation at $p$ gives a natural
homomorphism of sheaves $\co(-1) \to \cl|_{\Jac_d \times p} = \co$.
This vanishes precisely on $\Si_{d-1}$, so $\co(1) =
\co(\Si_{d-1})$.

To prove (ii), first assume $d>2g-2$.  Now for any $L \in \Jac_d$, the
restriction of $\co(\Delta)$ to $\Pj H^0(L) \times \Si$ is $L(1)$,
because $L(1)$ has a section with the universal property.  Hence
$\cl(1) \otimes \co(-\Delta)$ is trivial on the fibers of the
projection $\Si_d \times \Si \to \Jac_d$, which is locally trivial and
hence flat for $d>2g-2$.  Therefore $\cl(1) \otimes \co(-\Delta)$ is
the pull-back of some line bundle on $\Jac_d$; see Hartshorne
\cite[III Ex.\ 12.4]{h}.  But the restriction of $\co(\Delta)$ to
$\Si_d \times p$ is $\co(\Si_{d-1}) = \co(1)$ by (i), and the
restriction of $\cl$ to $\Si_d \times p$ is $\co$ by construction, so
this line bundle is trivial.

The case $d \leq 2g-2$ follows from this one by descending induction,
since the embedding $\iota: \Si_{d-1} \hookrightarrow \Si_d$ satisfies
$\iota^* \cl = \cl(p)$, $\iota^* \co(1) = \co(1)$, and $\iota^*
\co(\Delta) = \co(\Delta)(p)$.

Parts (iii) and (iv) then follow from (ii) together with well-known
formulas for $c_1(\cl)$ and $c_1(\co(\Theta))$; see \S\S2.6 and 2.7 of
Griffiths-Harris \cite{gh}.  \fp

A presentation of the cohomology ring of $\Si_d$ was given by
Macdonald \cite{mac}.

\bs{Theorem (Macdonald)}
\label{b}
The cohomology ring $H^*(\Si_d;\Z)$ is generated by $\xi_i$ and $\eta$
with relations
$$0 = \eta^r \prod_{i \in I} (\eta - \si_i) \prod_{j \in J} \xi_j
\prod_{k \in K} \xi_{k+g},$$
where $I,J,K \subset \{ 1, \dots, g \}$
are disjoint and
\begin{equation}
\label{a}
r + 2|I| + |J| + |K| \geq d+1.
\end{equation}
For $d > 2g-2$, these relations are generated by the single relation
\begin{equation}
\label{c}
0 = \eta^{d-2g+1}\prod_{i=1}^g (\eta - \si_i). 
\end{equation}
For $d \leq 2g-2$, they are generated
by those for which $r=0$ or $1$ and equality holds in \re{a}.
\fp \es

(Macdonald's paper contains a small error: it is asserted that $r=0$
is enough in the last statement.  However, the method of proof clearly
requires $r=1$ as well, and it is already necessary in the case
$d=2$.)

\bs{Corollary}
\label{gg}
The subring of $H^*(\Si_d)$ invariant under all monodromies of $C$
through smooth curves is generated by $\eta$ and $\theta$.
\es

\pf.  Certainly the monodromy invariant part of $H^*(\Jac_d)$ is
generated by $\theta$.  Indeed, since the monodromy surjects on $\Aut
H^*(\Si) = \Sp{2g,\Z}$, and $H^*(\Jac_d)$ is the exterior algebra on
$H^1(\Si)^*$, this means simply that the symplectic form and its powers
are the only alternating forms invariant under the symplectic group.
Macdonald's result shows that $H^*(\Si_d)$ is generated by $\eta$ as
an algebra over $H^*(\Jac_d)$.  For $d > 2g-2$, $\eta$ is monodromy
invariant by \re{x}(i); for $d<2g-2$, it is still monodromy invariant
since the embedding $i: \Si_{d-1} \hookrightarrow \Si_d$ satisfies
$i^* \eta = \eta$.  \fp

\smallskip 

We now turn to the quantum product on $\Si_d$.  First notice that
since $\si_i \in H^1(\Si_d)$, the quantum product with $\xi_i$ equals
the cup product by \re{d}(ii).  Hence the quantum product is completely
determined by the values of $\eta^u * \eta^v$ for $u,v \geq 0$.

Also notice that for $d>1$, $h_2(\Si_d) = {2g \choose 2} + 1$.  This
is an alarmingly large number of deformation parameters, but in fact
only one parameter is nontrivially involved in the quantum product,
for the following reason.

\bs{Proposition}
\label{z}
If $\ell \in H_2(\Si_d; \Z)$ is the homology class of any line in any
fiber $AJ^{-1}(L)$, then the quantum product on $H^*(\Si_d; \La)$
preserves the subring $H^*(\Si_d, \Q[[q^\ell]])$.  
\es

\pf.  Since an abelian variety has no rational curves whatsoever,
every genus 0 stable map to $\Si_d$ has image contained in a fiber of
$AJ$.  But for any line $\ell$ in any fiber of $AJ$, clearly $\eta
\cdot \ell = 1$, while $\nu \cdot \ell = 0$ for any class $\nu$ pulled
back from $\Jac_d$.  It follows that all such lines are homologous, so
every genus 0 stable map has image homologous to a non-negative multiple
of $\ell$.  \fp

\br{Remarks}

(i) Though it is not needed in the sequel, a topological version of
this statement remains true: the Hurewicz homomorphism $\pi_2(\Si_d)
\to H_2(\Si_d; \Z)$ has rank 1 for $d>1$, so the multiples of $\ell$
are the only spherical classes. This is proved in an appendix, \S9.

(ii) This subring determines the whole quantum product by
$\La$-linearity, so from now on attention will focus on it alone.  In
particular, the choice of polarization on $\Si_d$ used to define $\La$
is immaterial.  By abuse of notation $q^\ell$ is henceforth denoted
simply by $q$, $\langle \cdots \rangle_\ell$ by $\langle \cdots
\rangle_1$, and so on.

(iii) Macdonald \cite{mac} also shows that $c_1(T \Si_d) = (d-g+1)
\eta - \theta$; consequently, the degree of $q \in QH^*(\Si_d)$ is $2
c_1(T\Si_d) \cdot \ell = 2(d-g+1)$.  In particular, it is negative for
$d < g-1$ and 0 for $d = g-1$.  The latter case resembles that of a
Calabi-Yau manifold: although the canonical bundle is not trivial, its
restriction to every rational curve is trivial.  \er

\br{Example} 
\label{p}
Take $g=2$, $d=2$.  Then by Riemann-Roch, every line bundle in
$\Jac_2$ has one section, except the canonical bundle, which has two.
The Abel-Jacobi map therefore collapses exactly one rational curve $E
= \Pj H^0(K)$.  It is therefore precisely the blow-down of $E$ \cite[V
5.4]{h}.  The Poincar\'e dual of $E$ is easily seen to be $\theta -
\eta$.  Since $E$ is the only rational curve on $\Si_2$,
$\mbar_0(\Si_2,[E])$ is a point and $\mbar_{0,1}(\Si_2,[E]) = E$.
Therefore $\eta * \eta = \eta^2 + q (\theta - \eta)$, which completely
characterizes the quantum product. 
\er

In some cases the vanishing of the 3-point invariants, and hence of
certain terms in the quantum product, follows immediately from a
dimension count.

\bs{Proposition}
\label{hh}
For $d < g-1$ and $e > \frac{d-3}{g-1-d}$, $\langle a,b,c \rangle_e =
0$ for all $a,b,c \in H^*(\Si_d)$.
\es

\pf.  The expected dimension of $\mbar_0(\Si_d, e)$ is
$d-3+e(d-g+1)$, which is negative in this case.  The result follows
from \re{h}(iii).  \fp

Indeed, not only the 3-point invariants, but all higher-point
invariants vanish in this range by the same argument.

\bs{Corollary}
\label{bb}
For $d< g-1$ and $e > \frac{d-3}{g-1-d}$, $\Coeff_{q^e} a * b = 0$. In
particular, for $d < g/2 + 1$, the quantum product is simply the
ordinary cup product.  \fp
\es

It is also relatively easy to show that the 3-point invariants vanish
for $d > 2g-2$ and $e>1$.  Indeed, this follows from the vanishing of
the higher degree equivariant 3-point invariants of projective space.
However, we will not pursue this now, as it is subsumed in
\re{ii}(i) below.

\bit{The Brill-Noether approach}

To compute Gromov-Witten invariants for $\Si_d$ for $d\leq 2g-2$, we
must understand in detail how the fibers of the Abel-Jacobi map fit
together.  More precisely, we must understand the enumerative geometry
of the strata where the dimension of the fiber is constant.  This is
the subject matter of Brill-Noether theory.  Some well-known
definitions and results in the theory are recalled in \re{k}--\re{cc}
below.

We shall also make crucial use of a formula of Harris and Tu \cite{ht}
for the Chern numbers of kernel and cokernel bundles on determinantal
varieties.  Harris and Tu prove slightly more than the main result
they state; moreover, their statement contains a sign error (the very
first $+$ sign in the paper should be a $-$).  In the form we shall
need it, the result is the following.

\bs{Theorem (Harris-Tu)} 
\label{o}
Let $M$ be a complex manifold, $E$ and $F$ locally free sheaves of
ranks $m$ and $n$, and $f \in H^0(M, \itHom(E,F))$.  Let $S$ be the
universal subbundle over $\Gr_k E$, and let $\tilde{f} \in H^0(\Gr_k
E; \itHom(S,F))$ be the induced map.  Suppose $\tilde{f}$ intersects
the zero-section transversely in a variety $M_k$.  If $x_1, \dots,
x_k$ are the Chern roots of $S^*$ on $M_k$, then any characteristic
number $\nu c_1^{j_1}(S) c_2^{j_2}(S) \cdots c_k^{j_k}(S) [M_k]$,
where $\nu \in H^*(M)$, can be calculated using the formal identity
$$\nu x_1^{i_1}x_2^{i_2}\cdots x_k^{i_k} [M_k] = 
\nu \left| \begin{array}{llcl}
c_{n-m+k+i_1} & c_{n-m+k+1+i_1} & \cdots & c_{n-m+2k-1+i_1}\\
c_{n-m+k-1+i_2} & c_{n-m+k+i_2} & & \\
\phantom{xx} \vdots & & \ddots & \\
c_{n-m+1+i_k} & & & c_{n-m+k+i_k}
\end{array} \right| [M],$$
where $c_\al = c_\al(F-E)$. \fp
\es 

\br{Definitions} 
\label{k}
The moduli space $G^r_d$ of (possibly incomplete) linear systems on
$\Si$ of dimension $r$ and degree $d$ can be constructed just as in
the above theorem, taking $M = \Jac_d$.  Fix any reduced divisor $P$
on $\Si$ of large degree, and let $\cl$ be the Poincar\'e line bundle
on $\Jac_d \times \Si$.  Then the natural map $\cl(P) \to \co_P
\otimes \cl(P)$ pushes forward to a map $E \to F$ of locally free
sheaves on $\Jac_d$ such that, for any $r$, $G^r_d$ is precisely the
locus $M_k$ of the theorem for $k = r+1$.  The tautological subbundle
over $\Gr_{r+1} E$ restricts to a bundle $U^r_d$ over $G^r_d$ whose
projectivization $\Pj U^r_d$ admits a canonical map $\tau: \Pj U^r_d
\to \Si_d$.  Its image is referred to as $\Si^r_d$, and the subscheme
$AJ(\Si^r_d)$ of $\Jac_d$ is referred to as $W^r_d$.  Most of these
definitions are discussed at greater length by Arbarello et al.\ 
\cite{acgh}.  \er

\bs{Lemma}
\label{ff}
The total Chern classes of the bundles $E$ and $F$ defined above are
$c(E) = \exp(-\theta)$ and $c(F) = 1$.
\es

\pf. The Poincar\'e bundle is normalized so that for some $p \in \Si$,
$\cl|_{\Jac_d \times p} = \co$.  Since $F$ is just a sum of bundles
deformation equivalent to this, $c(F) = 1$.  As for $c(E)$, this can be
calculated using Grothendieck-Riemann-Roch: see Arbarello et al.
\cite[VIII \S2]{acgh}.  \fp

Since Gromov-Witten invariants are deformation invariants, the curve
$\Si$ may be taken to be general, and we will assume this henceforth.
This allows us to use one of the central results of Brill-Noether
theory.

\bs{Theorem (Gieseker)}
\label{s}
Let $\Si$ be a general curve.  For any effective divisor $D$ on $\Si$,
the natural map
$$H^0(\co(D)) \otimes H^0(K(-D)) \to H^0(K)$$
is injective. \es

\pf. See Gieseker \cite{gie}.  \fp

\bs{Corollary}
\label{cc} 
For a general curve $\Si$, $G^r_d$, and hence $W^r_d \sans W^{r+1}_d$
and $\Si^r_d \sans \Si^{r+1}_d$, are smooth of dimension $\rho$,
$\rho$, and $\rho+g$ respectively, where $\rho$ is the Brill-Noether
number $g-(r+1)(g-d+r)$. 
\es

\pf.  See Arbarello et al. \cite[V 1.6]{acgh}.  \fp

A further useful consequence of Gieseker's result is the following.

\bs{Lemma} 
\label{r}
\bl
\item At $L \in W^r_d \sans W^{r+1}_d$, there is a natural isomorphism
$N_{W^r_d/\Jac_d} = \Hom(H^0(L), H^1(L))$. 

\item At $D \in \Si^r_d \sans \Si^{r+1}_d$, there is a natural
  isomorphism 
$$ N_{\Si^r_d/\Si_d} = 
\Hom{\Big (}\mbox{\large $\frac{H^0(\co(D))}{\langle D \rangle}$},
H^1(D) {\Big )}.$$
\el
\es

\pf.  Part (i) follows easily from Proposition 4.2(i) in Chapter IV
of Arbarello et al. \cite{acgh}.  Part (ii) follows in the same way
from Lemma 1.5 in that chapter, provided that the composition of
natural maps
$$H^0(\co(D)) \otimes H^0(K(-D)) \to H^0(K) \to H^0(K \otimes \co_D)$$
has kernel $\langle D \rangle \otimes H^0(K(-D))$.  But the latter map
has kernel exactly $H^0(K(-D))$, so this follows from Gieseker's
theorem.  \fp

\bs{Proposition}
\label{f}
For general $\Si$, the reduced induced subscheme of the moduli space
$\mbar_0(\Si_d,e)$ is a disjoint union of orbifolds
$$\mbar_0(\Si_d,e)^\red 
= \bigcup_{r=1}^e \, \, \bigcup_{i=1}^r M_{r,i},$$
consisting of the stable maps whose image spans a linear system of
dimension $i$ contained in a complete linear system of dimension $r$
(or $\geq r$ if $r = e$).  The closure of each stratum has a
resolution $\tilde{M}_{r,i}$ which is the $\mbar_0 (\Pj^i,e)$ bundle
associated to the tautological subbundle over $\Gr_{i+1}U^r_d$.
\es

\pf.  The resolutions $\tilde{M}_{r,i}$ can be viewed as moduli spaces
of triples consisting of a linear system of dimension $r$, a
projective subspace
of dimension $i$, and a stable map to that subspace.  By 
\re{cc} they are orbifolds.  The open subset where the map spans the
$i$-dimensional subspace is a subvariety of $\mbar_0(\Si_d,e)$.  These
subvarieties partition $\mbar_0(\Si_d,e)$, because a rational curve of
degree $e$ in projective space spans at most an $e$-dimensional
subspace.  \fp

Since every complete linear system has dimension $\geq d-g$, $M_{r,i}
= \emptyset$ if $r<d-g$.  Otherwise $M_{r,i} \neq \emptyset$, and
$M_{r,i-1}$ is in its closure.  Hence the irreducible components of
$\mbar_0(\Si_d, e)^\red$ are the closures of $M_{r,r}$ for
$\min(d-g,e) \leq r \leq e$.

\bs{Proposition}
\label{g}
At any point of $M_{r,r}$, $\mbar_0(\Si_d, e)$ is reduced.
\es

Before the proof, a few observations that will be useful again later.

\bs{Lemma}
\label{q}
For any stable map $\phi: \curve \to \Si_d$, 
\bl
\item the natural map $T^1_{\Si_d}(\phi) \to T^1(\curve)$ of deformation
spaces is surjective; and

\item if $\Pj H^0(L)$ is the complete linear system containing
$\phi(\curve)$, the obstruction space $T^2_{\Si_d}(\phi)$ is
naturally isomorphic to $H^1(\curve,\phi^*N_{\Pj H^0(L)/\Si_d})$.
\el 
\es

\pf.  Any stable map has image in some $\Pj^r = \Pj H^0(L) \subset
\Si_d$, so there is a natural diagram of deformation spaces
$$\renewcommand{\arraystretch}{1.3} 
\begin{array}{ccccccccc} 
0 & \lrow & T^0(\curve) & \lrow & H^0(\curve, \phi^* T \Pj^r) & \lrow &
T^1_{\Pj^r}(\phi) & \lrow & T^1(\curve) \\
& & \bigdowneq & & \bigdownarg{} & & \bigdownarg{} & & \bigdowneq \\
0 & \lrow & T^0(\curve) & \lrow & H^0(\curve, \phi^* T \Si_d) & \lrow &
T^1_{\Si_d}(\phi) & \lrow & T^1(\curve). 
\end{array}
$$
By \re{j}, $H^1(\curve, \phi^* T\Pj^r) = 0$, so by \re{h}(i) the last
arrow in the first row is surjective, which proves (i).  Then
\re{h}(i) also implies that $T^2_{\Si_d}(\phi) =
H^1(\curve,\phi^*T\Si_d)$, and the long exact
sequence of
$$ 0 \lrow \phi^*T\Pj^r \lrow \phi^*T\Si_d 
\lrow \phi^*N_{\Pj^r/\Si_d} \lrow 0$$ 
completes the proof of (ii). \fp

\noindent{\it Proof of \re{g}}.  It suffices to show that for $\phi
\in M_{r,r}$, the deformation space $T^1_{\Si_d}(\phi)$ has dimension
equal to that of $M_{r,r}$ itself.  Now $M_{r,r}$ is an
$\mbar_0(\Pj^r, e)$ bundle over an open set in $G^r_d$, so
$$\dim M_{r,r} = \dim \mbar_0(\Pj^r, e) + \dim G^r_d.$$
On the other hand, the number of deformations of $\phi$ as a map to
$\Pj^r$ equals $\dim \mbar_0(\Pj^r,e)$ since the obstructions vanish.
By \re{q}(i), the number of deformations as a map to $\Si_d$ exceeds
this by 
$$\dim H^0(\curve,\phi^*T\Si_d) - \dim H^0(\curve,\phi^*T \Pj^r) 
= \dim H^0(\curve,\phi^*N_{\Pj^r/\Si_d}).$$
Hence it suffices to show $\dim H^0(\curve, \phi^* N_{\Pj^r/\Si_d}) =
\dim G^r_d$.  Because of the definition of $M_{r,r}$, there are two
cases, depending on whether or not $r=e$.  

If $r<e$, the image of $\phi$ is in $\Si^r_d \sans \Si^{r+1}_d$, which
is a $\Pj^r$-bundle over the open set $W^r_d \sans W^{r+1}_d \subset
G^r_d$.  On any fiber $\Pj^r = \Pj H^0(L)$ of this bundle, there is an
exact sequence
$$ 0 \lrow T_LW^r_d \otimes \co \lrow N_{\Pj^r/\Si_d} \lrow
N_{\Si^r_d/\Si_d} \lrow 0.$$
Hence it suffices to show that $H^0(\curve,\phi^* N_{\Si^r_d/\Si_d}) =0$.
By Lemma \re{r}(ii), there is a short exact sequence on $\Pj^r$
$$ 0 \lrow N_{\Si^r_d/\Si_d} \lrow \co \otimes \Hom (H^0(L), H^1(L))
\lrow \co(1) \otimes H^1(L) \lrow 0.$$
By the definition of $M_{r,r}$, the image of $\phi$ spans $\Pj
H^0(L)$.  Hence the natural map $H^0(L)^* \to H^0(\curve,
\phi^*\co(1))$ is injective, so 
$$H^0(\curve, \co \otimes \Hom (H^0(L), H^1(L))) 
\lrow H^0(\curve, \phi^* \co (1) \otimes H^1(L))$$ 
is also and $H^0(\curve, \phi^* N_{\Si^r_d/\Si_d}) = 0$ as desired.  

If $r=e$, then the image of $\phi$ spans a linear system $\Pj^r$ which
may not be complete.  Let $\Pj H^0(L)$ be the complete linear system 
containing it, and let $\ell$ be its projective dimension.  Then there
are three short exact sequences on $\Pj^r$:
$$0 \lrow N_{\Pj^r/\Pj H^0(L)} \lrow N_{\Pj^r/\Si_d} 
\lrow N_{\Pj H^0(L)/ \Si_d} \lrow 0, $$
$$0 \lrow T_LW^\ell_d \otimes \co \lrow N_{\Pj H^0(L)/ \Si_d} 
\lrow N_{\Si^\ell_d/\Si_d} \lrow 0, $$
$$0 \lrow N_{\Si^\ell_d/\Si_d} \lrow \Hom(H^0(L),H^1(L)) \otimes
\co \lrow H^1(L) \otimes \co(1) \lrow 0,$$
the last by \re{r}(ii).  But $N_{\Pj^r/\Pj H^0(L)} \cong \co^{\ell -
r}(1)$ since $\Pj^r \subset \Pj H^0(L)$ is a projective subspace.
Since $H^1(\curve, \co) = H^1(\curve, \phi^*\co(1)) = 0$, all three of
the long exact sequences contain only $H^0$ terms.  Hence
$\dim H^0(\curve, \phi^* N_{\Pj^r/\Si_d})$ can be expressed in terms
of $\dim H^0(\curve, \phi^* \co(1))$, which is $r+1$ since $\phi$
spans $\Pj^r$, and the ranks of the bundles above.  A little
high-school algebra gives the desired result. \fp

\bit{The degree 1 invariants}

\bs{Proposition}
\label{ee}
As schemes, $\mbar_0(\Si_d, 1) = G^1_d$ and $\mbar_{0,1}(\Si_d,1)
= \Pj U^1_d$.
\es

\pf.  The identification of sets is clear from \re{f}, and  the moduli
space is reduced by \re{g}.  \fp

\bs{Theorem}
\label{l}
For $u,v,w \geq 0$,
$$\langle \eta^u, \eta^v, \theta^{g-m} \eta^w \rangle_1 
= \frac{g!}{m!} \sum_{i=0}^{u-1} 
\bino{m}{g-d+i+v}  
- \bino{m}{g-d+i},$$
where $m = 2g-2d-1+u+v+w$.
\es

For example, consider the case when the expected dimension $\rho = g -
2(g-d+1)$ of $G^1_d$ is zero, which occurs when $g=2d-2$.  Then
$\langle \eta, \eta, \eta \rangle_1$ counts the number of points in
$G^1_d$.  According to \re{l}, this is ${2d - 2 \choose d-1} - {2d-2
\choose d-2}$, which agrees with the Catalan number
$\frac{1}{d}{2d-2 \choose d-1}$ computed by Castelnuovo \cite[V
1.2]{acgh}.

\pf. Since the moduli space is smooth of the expected dimension, the
Gromov-Witten invariants can be calculated using \re{ll}.  The
evaluation map is precisely $\tau: \Pj U^1_d \to \Si_d$, as defined in
\re{k}.  By \re{x}(i) $\tau^* \eta = c_1(\co(1))$, and hence $f_* \tau^*
\eta^u$ is the Segre class $\seg_{u-1}(U^1_d)$; see for example Fulton
\cite[\S3.1]{f}.  If $x_1$ and $x_2$ are the Chern roots of
$(U^1_d)^*$, then the Segre classes are the complete symmetric
polynomials:
$$\seg_k = \sum_{i=0}^k x_1^i x_2^{k-i}.$$
So the 3-point invariant is
$$\langle \eta^u, \eta^v, \theta^{g-m} \eta^w \rangle_1 =
\theta^{g-m} \sum_{i=0}^{u-1} \sum_{j=0}^{v-1}
 \sum_{k=0}^{w-1} x_1^{i+j+k}
x_2^{u+v+w-3-i-j-k} [G^1_d].$$
The Harris-Tu formula implies that
$$\theta^{g-m}x_1^p x_2^{u+v+w-3-p} [G^1_d]
= \theta^{g-m}(\ga_p-\ga_{p+1}) [\Jac_d],$$
where 
\begin{eqnarray*}
\ga_p & = & c_{g-d+1+p}(F-E) \cdot c_{g-d+u+v+w-2-p}(F-E) \\
 & = &
\theta^m /(g-d+1+p)! (g-d+u+v+w-2-p)!,
\end{eqnarray*}
the last equality by \re{ff}.
The sum over $k$ therefore telescopes to yield
$$\langle \eta^u, \eta^v, \theta^{g-m} \eta^w \rangle_1 =
\theta^{g-m} \sum_{i=0}^{u-1} \sum_{j=0}^{v-1}
(\ga_{i+j} - \ga_{i+j+w}) [\Jac_d].$$ 
But there is another symmetry, namely $\ga_p = \ga_{u+v+w-3-p}$;
applying this to the second term and canceling gives
\begin{eqnarray*}
\sum_{i=0}^{u-1} \sum_{j=0}^{v-1}
\ga_{i+j} - \ga_{i+j+w} 
& = &
\sum_{i=0}^{u-1} \sum_{j=0}^{v-1}
\ga_{i+j} - \ga_{(u-1-i)+(v-1-j)-1} \\
& = &
\sum_{i=0}^{u-1} \sum_{j=0}^{v-1}
\ga_{i+j} - \ga_{i+j-1} \\
& = &
\sum_{i=0}^{u-1}
\ga_{i+v-1} - \ga_{i-1}.
\end{eqnarray*}

Substituting this and using $\theta^g [\Jac_d] = g!$ then gives the
answer as stated.  \fp

Of course, by the same method one could easily derive a formula for
3-point invariants involving arbitrary elements of $H^*(\Si_d)$.  But
these are enough to characterize the linear term of the quantum
product.

\bs{Corollary}
\label{t}
For $u,v \geq 0$, 
$$\eta^u * \eta^v = \eta^{u+v} 
+ q \sum_{i=0}^{u-1} \left( \frac{\theta^{g-d+i+v}}{(g-d+i+v)!}
\, \eta^{u-1-i} - \frac{\theta^{g-d+i}}{(g-d+i)!} \, \eta^{v+u-1-i}
\right) 
+ O(q^2).$$
\es

\pf.  By the definition of quantum product, the coefficient of $q$ in
$\eta^u * \eta^v$ is $\mu$ if and only if for all $\nu \in
H^*(\Si_d)$, $\langle \eta^u, \eta^v, \nu \rangle_1 = \mu \nu
[\Si_d]$.  However, Gromov-Witten invariants, and hence the quantum
product, are deformation invariant.  Since $\eta$ and its powers are
monodromy invariant by \re{gg} whenever $\Si$ lies in a family of
smooth curves, each coefficient of $\eta^u * \eta^v$ must be as well.
Again by \re{gg}, this means that it is in the subring generated by
$\eta$ and $\theta$.  This subring is the algebraic part of
$H^*(\Si_d)$ for $\Si$ general: see Arbarello et al.\ \cite[VIII
\S5]{acgh}.  In particular, it satisfies Poincar\'e duality; hence it
suffices to take $\nu$ in the subring as well.  This circumvents the
cumbersome task of dealing with arbitrary monomials in the $\xi_i$.

It is easy to check against \re{l} that the coefficient stated above
satisfies the required condition.  One simply applies the convenient
formula 
$$\theta^i \eta^{d-i} [\Si_d] = 
\left\{ \begin{array}{cl}
\frac{g!}{(g-i)!} & \mbox{if $i \leq g$}\\
0 & \mbox{if $i>g$}
\end{array} \right. $$
which follows from \re{x}(i) and \re{c} using the fact that
$\Si_{d-1} \subset \Si_d$ is Poincar\'e dual to $\eta$.  The required
identity follows term-by-term from a comparison of the coefficient
above with the Gromov-Witten invariant in \re{l}.  \fp

\bit{The degree 2 invariants}

The degree 2 invariants can also be calculated using ideas from
Brill-Noether theory.  Here things are considerably more complicated;
in particular, the virtual class comes into play.  However, some
delightful cancellations make the computations tractable.  

In the decomposition of \re{f}, the moduli space $\mbar_0(\Si_d,2)$ has
three strata, $M_{11}$, $M_{22}$, and $M_{12}$.  If $d\geq g+2$,
$M_{11}$ is empty, but this will not affect the results.

Of these strata, $M_{22}$ has the expected dimension, while $M_{11}$
exceeds it by $g+1-d$.  The third, $M_{12}$, is in the closure of
$M_{22}$ and has less than the expected dimension.  The virtual class
is therefore a sum of cycles pushed forward from $\tilde{M}_{11}$ and
$\tilde{M}_{22}$.  By \re{g}, the former is simply the ordinary
fundamental class of $\overline{M}_{22}$.  The Gromov-Witten invariant
therefore can be computed as in \re{ll}:
$$\langle a_1,a_2,a_3 \rangle_2 = 
{\Big (} \prod_k \tilde{f}_* \widetilde{\ev}^* a_k {\Big )}
([\tilde{M}_{11}]^\vir + [\tilde{M}_{22}]).$$ 

\bs{Proposition}
\label{m}
For all $a_1,a_2,a_3 \in H^*(\Si_d)$, 
$${\Big (} \prod_k \tilde{f}_* \widetilde{\ev}^* a_k 
{\Big )}[\tilde{M}_{22}] = 0. $$
\es

\pf.  On $\tilde{M}_{22} \times_{\mbar_0(\Si_d,2)}
\mbar_{0,1}(\Si_d,2)$, the evaluation maps factor through the $\Pj^2$
bundle $\Pj U^2_d$, where any cohomology class can be expressed as
a polynomial over $H^*(G^2_d)$ of degree $\leq 2$.  So on
$\tilde{M}_{22}$, $\tilde{f}_* \, \widetilde{\ev}^* a_i$ is in the
submodule $H^{\geq \deg a_i - 4}(G^2_d) \cdot H^*(\tilde{M}_{22})$.
It certainly vanishes on the fundamental class for dimensional reasons
unless $\sum \deg a_i = 2 \dim \tilde{M}_{22} + 6 = 6d-4g+4$.  But
then the product is in the submodule $H^{\geq 6d-4g-8}(G^2_d) \cdot
H^*(\tilde{M}_{22})$, which is zero since $\dim G^2_d$ is the
Brill-Noether number $\rho = 3d-2g-6$. \fp

The above is really a special case of the more general statement that
the equivariant 3-point invariants of projective space vanish in
degree $>1$.

We now attack the virtual class $[\tilde{M}_{11}]^\vir$, or rather its
push-forward to $G^1_d$, where the invariant will be calculated.  For
simplicity $U^1_d$ will be denoted $U$ in the rest of this section.

\bs{Theorem} 
\label{n}
Let $p: \tilde{M}_{11} \to G^1_d$ be the projection.  
Then $p_* [\tilde{M}_{11}]^\vir$ is Poincar\'e dual to the degree
$2(g-1-d)$ part of
$$\frac{\ratio{1}{8} \exp \theta}{1+ c_1(U)/2}.$$
\es

\pf.  First, it suffices to work only on $p: M_{11} \to W^1_d
\sans W^2_d$, because the codimension of the missing locus in $G^1_d$
is $g+2-d$, which is greater than $g-1-d$.

Take the short exact sequence
$$0 \lrow \co \lrow \co(\Delta) \lrow \co_\Delta(\Delta) \lrow 0$$ on
$(\Si^1_d \backslash \Si^2_d) \times \Si$, and push forward to a
long exact sequence on the first factor:

$$ \begin{array}{ccccccl} 
0 & \lrow & (R^0 \pi)_* \co & \lrow & (R^0 \pi)_*
\co(\Delta) & \lrow & (R^0 \pi)_* \co_\Delta (\Delta) \\  
 & \lrow & (R^1 \pi)_* \co & \lrow & (R^1 \pi)_*
\co(\Delta) & \lrow & 0.
\end{array} $$
The third nonzero term is $T\Si_d$, and the image of the second is the
tangent space to the pencils.  The quotient is therefore the vector
bundle $N$ whose fibers along each pencil are the normal spaces to
that pencil in $\Si_d$.  So there
is a short exact sequence on $\Si^1_d \sans \Si^2_d$
$$0 \lrow N \lrow (R^1 \pi)_* \co \lrow (R^1 \pi)_* \co(\Delta) \lrow
0,$$
and hence a long exact sequence
$$0 \lrow H^0(\curve, \phi^* N) \lrow H^1(\Si, \co) \lrow
H^0(\curve, \phi^* (R^1 \pi)_* \co(\Delta)) 
\lrow H^1(\curve, \phi^* N) \lrow 0.$$
By \re{q}(ii), the obstruction space $T^2_{\Si_d}(\phi)$ is nothing
but the last term of this sequence.  Moreover, by \re{h}(ii), the
virtual class on $M_{11}$ is exactly the Euler class of the orbifold
vector bundle whose fiber at $\phi: \curve \to \Si_d$ is this
obstruction space.

How can the terms in this long exact sequence be described better?
Well, note that $N$ also fits into the exact sequence on $\Si^1_d
\backslash \Si^2_d$
$$0 \lrow TW^1_d \lrow N \lrow N_{\Si^1_d/\Si_d} \lrow 0.$$
By \re{r}(ii), $N_{\Si^1_d/\Si_d}$ restricted to the pencil $\Pj H^0(L)$
is $\co(-1) \otimes H^1(L)$, so $H^0(\curve,\phi^* N_{\Si^1_d/\Si_d})
$ vanishes and hence $H^0(\curve, \phi^* N)$ is simply $T_LW^1_d$.  On
the other hand, if $\co(1)$ on $\Si_d$ denotes the relative hyperplane
bundle of the Abel-Jacobi map, and $\cl$ is a Poincar\'e line bundle
pulled back from $\Jac_d \times \Si$, then on $\Si_d \times \Si$,
$\co(\Delta) = \cl(1)$ by \re{x}(ii).  Hence on $\Si^1_d \sans \Si^2_d$,
$(R^1\pi)_* \co(\Delta)$ is the tensor product $V(1)$, where $V = (R^1
\pi)_* \cl$ is a vector bundle pulled back from $W^1_d \sans W^2_d$.
The term $H^0(\phi^* (R^1 \pi)_* \co(\Delta))$ which appeared above is
therefore the tensor product $V \otimes H^0(\curve, \phi^* \co(1))$.

Putting it all together yields an exact sequence
$$0 \lrow N_{W^1_g/\Jac_d} \lrow V \otimes H^0(\curve, \phi^* \co(1))
\lrow H^1(\curve, \phi^* N) \lrow 0$$
for any $\phi \in M_{11}$.  The first two terms clearly have dimension
independent of $\phi$, so they determine an orbifold vector bundle on
$M_{11}$, and hence so does the third term.  

Now consider the orbifold structure on $M_{11}$.  It is an
$\mbar_0(\Pj^1,2)$-bundle over $W^1_d \sans W^2_d$.  In fact, a degree
2 stable map to a line is characterized by its base points, so as
schemes $\mbar_0(\Pj^1,2) \cong \Pj^2$ and $M_{11} = \Pj \Sym^2 U^*$
where $U = U^2_d$.  But every stable map in $M_{11}$ has an
involution, so as an orbifold $M_{11}$ is a global quotient by $\Z_2$
of a double cover, say $\hat M_{11}$, branched over the bundle of
conics $\Pj U^* \subset \Pj \Sym^2 U^*$.

The involution splits $H^0(\curve, \phi^* \co(1))$ into $\pm
1$-eigenspaces.  The $+1$-eigenspace consists of those sections pulled
back from the pencil, so the corresponding bundle is just $U^*$.  The
$-1$-eigenspace is generated by the square root of the section of
$\co(2)$ on the pencil vanishing on the branch points of the double
cover, so the corresponding bundle is an orbifold line bundle on
$M_{11}$, coming from a $\Z_2$-equivariant line bundle on $\hat
M_{11}$ whose tensor square is the pull-back of $\co(-1)$ from
$M_{11}$.  Call this orbifold bundle 
$\co(-\half)$.

Now it follows from \re{r}(i) that the normal bundle
$N_{W^1_g/\Jac_d}$ is $V \otimes U^*$, and its image in $V \otimes
H^0(\curve, \phi^* \co(1))$ is clearly in the $+1$-eigenspace.  This
splits the exact sequence mentioned above.  Hence the orbifold bundle
with fiber $H^1(\curve, \phi^* N)$ is none other than $V(-\half)$, and
the virtual class is the Euler class of this.

It only remains to push forward the virtual class to $G^1_d$.  If
$\xi_j$ are the Chern roots of $V$, and $h$ is the hyperplane class of
$\Pj \Sym^2 U^*$, then the virtual class is $\prod_j (\xi_j - h/2)$.
This pushes forward to 
$$\half \sum_{n=0}^{g-d-1} 
(-\half)^{n + 2} c_{g-1-d-n}(V) \seg_n(\Sym^2 U^*),$$
where $\seg_i$ denotes the Segre class.  The extra factor of $\half$
appears because of the orbifold structure on $M_{11}$.

The Chern roots of $\Sym^2 U^*$ are $2x_1$, $2x_2$, and $x_1+x_2$.  The
factors of 2 cancel and one gets the degree $2(g-d-1)$ part of:
$$ \frac{\ratio{1}{8} \, c(V)}{(1-x_1)(1-x_2)(1-(x_1+x_2)/2)}.$$

Now a miraculous cancellation.  The $(1-x_1)(1-x_2)$ in the
denominator is just $c(U)$.  But $U$ and $V$ are exactly the kernel
and cokernel of the map of bundles $E \to F$ defined in \re{k}.  Hence
by \re{ff}
\begin{eqnarray*}
c(V)/c(U) & = & c(F) / c(E) \\
& = & 1/\exp(-\theta)\\
& = & \exp \theta,
\end{eqnarray*}
and plugging this in yields the stated formula.  \fp

\bs{Theorem}
For $u,v,w \geq 0$, 
\begin{eqnarray*}
\lefteqn{\langle \eta^u, \eta^v, \theta^{d+1-m} \eta^w \rangle_2} \\
& & \displaystyle = \sum_{n=0}^{g-1-d} \sum_{p=0}^n 
\bino{n}{p}\frac{g!}{2^n (g-1-d-n)! (m+n)!}
\sum_{i=0}^{u-1}
\bino{m + n}{g-d+i+v+p} - \bino{m + n}{g-d+i+p},
\end{eqnarray*}
where $m = 2g-2d-1+u+v+w$.
\es

\pf.
The proof is parallel to that of \re{l}.

By \re{m}, 
$$ \langle a_1, a_2, a_3 \rangle_2 = 
{\Big (} \prod_k \tilde{f}_* \widetilde{\ev}^* a_k 
{\Big )} [\tilde{M}_{11}]^\vir.$$

Conveniently, $\tilde{M}_{11} = \Pj \Sym^2 U^*$ is isomorphic to the
closure $\overline{M}_{11}$.  Also, $f^{-1}(\overline{M}_{11})$ is
simply the fibered product $\Pj \Sym^2 U^* \times_{G^1_d} \Pj U$.  In
contrast to $M_{11}$, a generic stable map in
$f^{-1}(\overline{M}_{11})$ has no involution; it is rigidified by the
marked point.  The evaluation map $f^{-1}(\overline{M}_{11}) \to
\Si_d$ factors through $\Pj U$; hence $f_* \ev^* \eta^u$ restricted to
$\overline{M}_{11}$ is a class pulled back from $G^1_d$, namely $2
\seg_{u-1}(U)$, where $\seg_i$ is the $i$th Segre class. The factor of
$2$ appears because of the orbifold structure on $\overline{M}_{11}$.
In terms of the Chern roots $x_1, x_2$ of $U^*$, this is
$$2 \seg_{u-1}(U) = 2 \sum_{i=0}^{u-1} x_1^i x_2^{u-1-i}.$$ 

Expanding the formula of \re{n} in terms of $x_1$ and $x_2$ yields the
cap product
$$p_* [\tilde{M}_{11}]^\vir = \ratio{1}{8} \sum_{n =
0}^{g-1-d} \sum_{p = 0}^n \bino{n}{p} \frac{(-1)^n
\theta^{g-1-d-n} x_1^p x_2^{n-p}}{2^n (g-1-d-n)!} \cap [G^1_d].$$

Consequently
\begin{eqnarray*}
\lefteqn{\langle \eta^u, \eta^v, \theta^{g-n-m} \eta^w \rangle_2} \\
& & \displaystyle = 
\sum_{n = 0}^{g-1-d} \sum_{p=0}^n
\sum_{i=0}^{u-1} \sum_{j=0}^{v-1} \sum_{k=0}^{w-1} 
\bino{n}{p}
\frac{\theta^{d-4-m} x_1^{p + i + j + k}x_2^{n -
p + u + v +w - 3 - i - j - k}}{2^n (g-1-d-n)!}[G^1_d].
\end{eqnarray*}

As in the proof of \re{l}, one now applies the Harris-Tu formula \re{o},
telescopes the sum over $k$, and cancels using the additional symmetry
to obtain the stated result.  \fp

\bs{Corollary}
\label{dd}
For $u,v \geq 0$, 
\begin{eqnarray*}\lefteqn{\eta^u * \eta^v = \eta^{u+v} 
+ q \sum_{i=0}^{u-1} \left( \frac{\theta^{g-d+i+v}}{(g-d+i+v)!}
\, \eta^{u-1-i} - \frac{\theta^{g-d+i}}{(g-d+i)!} \, \eta^{v+u-1-i}
\right)} \\
& & \displaystyle + \, \, q^2 \sum_{n=0}^{g-1-d} \sum_{p=0}^n 
\bino{n}{p}\frac{1}{2^n (g-1-d-n)!} \\
& & \displaystyle \cdot \sum_{i=0}^{u-1}
\left( \frac{\theta^{2g-2d+1-n+i+v+p}\eta^{u+n-i-p-3}}{(g-d+i+v+p)!} - 
\frac{\theta^{2g-2d+1-n+i+p}\eta^{u+v+n-i-p-3}}{(g-d+i+p)!} \right)
+ O(q^3).
\end{eqnarray*}
\es

\pf.  Similar to that of \re{t}.  \fp

\bit{The higher degree invariants for \boldmath $d \geq g-1$}

Examining the formulas of the last two sections reveals that the
degree 2 invariants vanish for $d > g-1$, and equal the degree 1
invariants for $d = g-1$.  The result below shows that this is the case
for all higher degree invariants.  

\bs{Theorem}
\label{ii} 
For all $e>1$ and all $a_1, a_2, a_3 \in H^*(\Si_d)$,
\bl
\item $\langle a_1, a_2, a_3 \rangle_e = 0$ when $d > g-1$;

\item $\langle a_1, a_2, a_3 \rangle_e 
= \langle a_1, a_2, a_3 \rangle_1$ when $d = g-1$. 
\el
\es

\pf. The virtual class is some algebraic cycle class of degree $d -3+
e(d-g+1)$.  Choose a cycle representing this class; it decomposes
uniquely as a sum of cycles $\sum_{r,i} [M_{r,i}]^\vir$, where every
component of $[M_{r,i}]^\vir$ is supported on the closure of $M_{r,i}$
but not on any smaller stratum.  Then $[M_{r,i}]^\vir$ is in the image
of the push-forward from the resolution $\tilde{M}_{r,i}$, so its
contribution to the Gromov-Witten invariant $\langle a_1,a_2,a_3
\rangle_e$ can be computed on $\tilde{M}_{r,i}$ as
$${\Big (} \prod_k \tilde{f}_*\widetilde{\ev}^* a_k 
{\Big )} [\tilde{M}_{r,i}]^\vir.$$

Now on $\tilde{M}_{r,i}$, all three of the evaluation maps factor
through $\Pj U^r_d$.  The pullback of $a_i$ to this bundle can be
written as a polynomial in the hyperplane class over $H^*(G^r_d)$ of
degree $\leq r$.  Hence $\tilde{f}_*\widetilde{\ev}^* a_i$ is in the
submodule
$$H^{\geq \deg a_i - 2r}(G^r_d) \cdot H^*(\overline{M}_{r,i}).$$
Certainly $\langle a_1, a_2, a_3 \rangle_e =0$ on dimensional grounds
unless $\sum_k \deg a_k = 2d + 2e(d-g+1)$.  In that case the product
$\prod_k \tilde{f}_*\widetilde{\ev}^* a_k$ is in the submodule
$$H^{\geq 2d + 2e(d-g+1) - 6r}(G^r_d) \cdot H^*(\overline{M}_{r,i}).$$
But the dimension of $G^r_d$ is the Brill-Noether number $\rho = g -
(r+1)(g-d+r)$, and some high school algebra shows that 
$$d + e(d-g+1) -3r - \rho = (e-r)(d-g+1) + r(r-1).$$  
This is positive for all $d > g-1$ unless $e=r=1$, and also for $d =
g-1$ unless $r=1$. 

Hence the ideal $H^{\geq 2d + 2e(d-g+1) - 6r}(G^r_d)$ vanishes in
those cases.  This immediately implies part (i) of the theorem, and it
shows that when $d=g-1$, only the single stratum $M_{11}$ contributes
nontrivially to the 3-point invariants.

We therefore turn to the contribution of this stratum, which consists
of $e$-fold covers of pencils on $\Si_{g-1}$. The closure
$\overline{M}_{11}$ is isomorphic to $\tilde{M}_{11}$, which is a
bundle over $G^1_d$ with fiber $\mbar_0 (\Pj^1,e)$.  The evaluation
map factors through $\Pj U^1_d$:
$$\begin{array}{rcccc}
f^{-1}(\overline{M}_{11}) & \lrow & \Pj U^1_d &
\stackrel{\tau}{\lrow} & \Si_d \\
 \leftbigdownarg{f}& & \bigdownarg{\pi} & & \\
 \overline{M}_{11}&\stackrel{p}{\lrow} & G^1_d. & &
\end{array}$$

The morphism $f^{-1}(\overline{M}_{11}) \to \overline{M}_{11}
\times_{G^1_d} \Pj U^1_d$ is finite of degree $e$, so for $a \in
H^*(\Si_d)$,
$$p_* f_* \ev^* a = e \cdot \pi_* \tau^* a.$$
The Gromov-Witten invariant can be expressed as
$$ \langle a_1, a_2, a_3 \rangle_e 
= {\Big (} \prod_k p_* f_* \ev^* a_k {\Big )} \, 
p_* [\tilde{M}_{11}]^\vir; $$ 
but since $G^1_d$ has dimension $g-4$, exactly that of the virtual
cycle, $p_* [\tilde{M}_{11}]^\vir \in H_{2g-8}(G^1_d)$ must be a
scalar multiple of $[G^1_d]$.

On the other hand, by \re{ee} $G^1_d = \mbar_0 (\Si_d,1)$ and $\Pj
U^1_d = \mbar_0 (\Si_d,1)$.  Hence
$${\Big (} \prod_k \pi_* \tau^* a_k {\Big )} [G^1_d] 
= \langle a_1, a_2, a_3 \rangle_1. $$ 
So to prove part (ii) of the theorem, it suffices to show that $p_*
[\tilde{M}_{11}]^\vir = 1/e^3 [G^1_d]$. One may work on any single
fiber of $p$; choose one over $W^1_d$.  The moduli space is an orbifold
in a neighborhood of this fiber, so \re{h}(ii) applies.  Hence we need
to contribute the Euler class of the orbifold bundle whose fiber at a map
$\phi$ is the obstruction space $T^2_{\Si_d}(\phi)$.  By \re{q}(ii),
$T^2_{\Si_d}(\phi) = H^1(\curve,\phi^* N_{\Pj^1/\Si_d})$.
From the short exact sequence
$$ 0 \lrow N_{\Pj^1/\Si^1_d} \lrow N_{\Pj^1/\Si_d} \lrow
N_{\Si^1_d/\Si_d} \lrow 0,$$
together with the isomorphism
$N_{\Pj^1/\Si^1_d} \cong \co^\rho$, it follows that 
$$H^1(\curve, \phi^* N_{\Pj^1/\Si_d}) 
= H^1(\curve, \phi^* N_{\Si^1_d/\Si_d}).$$  
As seen in  \re{r}(ii), the normal space to $\Si^r_d$ at a divisor
$D$ is naturally isomorphic to $\Hom(H^0(\co(D))/\langle D \rangle,
H^1(\co(D)))$.  In this case $\dim H^0 = \dim H^1 = 2$, so on $\Pj^1$
the normal bundle $N_{\Si^1_d/\Si_d}$ is isomorphic to $\Hom (\co(1),
\co \oplus \co) = \co(-1) \oplus \co(-1)$.

So, on $\mbar_0(\Pj^1,e)$, we want to know the Euler class of the
orbifold bundle whose fiber at a map $\phi$ is $H^1(\phi^* \co(-1)
\oplus \phi^*\co(-1))$.  Very felicitously, this is exactly the number
computed to be $1/e^3$, using some nontrivial combinatorics, in the
work of Aspinwall-Morrison \cite{am}, Manin \cite{m}, and Voisin
\cite{v}.  In their work, the motivation was to compute Gromov-Witten
invariants for Calabi-Yau threefolds containing rational curves with
normal bundle $\co(-1) \oplus \co(-1)$; the present case is in a sense
a relative version of this, since $\Si_{g-1}$ contains the family
$\Si^1_{g-1}$ of rational curves whose normal bundle restricts to
$\co(-1) \oplus \co(-1)$ on every curve.  In any case, this completes
the proof.  \fp

\bs{Corollary}
\label{aa}
For all $u,v \geq 0$, 
\bl
\item when $d > g-1$, 
$$\eta^u * \eta^v = \eta^{u+v}  
+ q \sum_{i=0}^{u-1} \left( \frac{\theta^{g-d+i+v}}{(g-d+i+v)!}
\, \eta^{u-1-i} - \frac{\theta^{g-d+i}}{(g-d+i)!} \, \eta^{u+v-1-i}
\right) ;$$

\item when $d = g-1$, 
$$\eta^u * \eta^v = \eta^{u+v} + \frac{q}{1-q}
\sum_{i=1}^u \left( \frac{\theta^{i+v}}{(i+v)!}
\, \eta^{u-i} - \frac{\theta^i}{i!} \, \eta^{v+u-i} \right) .$$
\el
\es

\pf.  Similar to that of \re{t}. \fp

\bit{Presentation of the quantum ring}

The results of the previous sections \re{gg}, \re{t}, \re{dd}, \re{aa}
have given an explicit quantum multiplication table for $\Si_d$
provided $d \not\in [\ratio{3}{4}g, g-1)$.  It is natural to look for
a presentation of the quantum cohomology ring as well.  Such a
presentation contains less information than the multiplication table,
because it does not specify the additive isomorphism between
$QH^*(\Si_d)$ and $H^*(\Si_d; \La)$.  We will be able to determine it
even in the slightly more general case $d \not\in
[\ratio{4}{5}g-\ratio{3}{5},g-1)$.  To go this far requires only the
degree 1 invariants and the results of the last section; in principle
one could continue as far as $\ratio{6}{7}g-\ratio{5}{7}$ by using the
degree 2 invariants, but this becomes cumbersome.

\bs{Proposition} 
The rings $QH^*(\Si_d)$ are generated over $\La$ by Macdonald's
generators $\eta$ and $\xi_i$, and there is a complete set of
relations, uniquely determined by the property that it reduces mod $q$
to Macdonald's relations.  
\es

\pf.  These facts are well-known for quantum cohomology generally
in the case $\deg q > 0$.  The first statement is proved by a simple
induction on the degree of the cup product of an arbitrary collection
of classical generators.  Just take a similar monomial where all the
cup products are replaced by quantum products. The difference between
the two monomials is a multiple of $q$, so the coefficient has lower
degree and by induction it can be expressed as a quantum product of
the classical generators.  This shows that the classical generators
generate the quantum ring as well.

To extend the classical relations to quantum relations, first take a
classical relation and replace all the cup products by quantum
products.  This quantum expression may not be a quantum relation, but
it is in the ideal $\langle q \rangle$, because it reduces to a
classical relation.  So express it as $q$ times a classical monomial,
replace the cup products in this monomial by quantum products, and
subtract the result from the quantum expression.  This difference is
now in $\langle q^2 \rangle$.  Proceed inductively; when $\deg q >0$,
the coefficient of a high power of $q$ will eventually be in $H^{<0} =
0$.

To see that these generate all the quantum relations, just notice that
imposing them gives a free module over $\La$ of the correct dimension.

Similar arguments work for $\deg q < 0$.  To prove the first
statement, for example, express the generator of top-dimensional
cohomology as a cup product of generators.  This identity still holds
if the cup products are replaced by quantum products.  Then apply
descending induction.

This leaves only the case $\deg q = 0$, which for a symmetric product
$\Si_d$ means $d=g-1$.  Since quantum and cup products with each
$\xi_i$ are the same, to prove the first statement it suffices to
express $\eta^u$ as a quantum polynomial in $\eta$ and $\xi_i$ for all
$u \geq 0$.  An induction on $u$ shows this is possible, since by
\re{aa}(ii)
$$\eta * \eta^u = \eta^{u+1} + \frac{q}{1-q} \left(
\frac{\theta^{u+1}}{(u+1)!} - \theta \eta^u \right).$$ 

The second statement is proved as for $\deg q >0$, but by induction on
the power of $\theta$ rather than $q$.  \fp

How do we find the quantum relations alluded to in the proposition?
Macdonald's classical relations from \re{b} can be expressed as 
\begin{equation}
\label{v}
0 = {\Big (} \sum_\al \eta^{r + |I| - \al} s_\al(\si_i) {\Big )}
\left( \prod \xi_j \prod \xi_{k+g} \right) 
\end{equation}
where $s_\al$ is the elementary symmetric polynomial of degree $\al$.
To find explicit quantum relations reducing to these mod $q$, the
following result is useful.  We adopt the notation $\prod^*$ for a
quantum product, and $\eta^{*u}$ for a quantum power.

\bs{Proposition}
\label{w}
For $u \geq 0$, $\eta^u =$
$$\setlength\arraycolsep{2pt}
\renewcommand{\arraystretch}{2.4}
\begin{array}{rcrclcl}
\rm{(i)} & &
\eta^{*u} & - & \displaystyle  q \sum_{j=0}^{n}
\eta^{*(n-j)}\theta^j/j! & &
\mbox{for $d > g$, where $n=u-d+g-1$;} \\
\rm{(ii)} & &
(\eta+q)^{*u} & - & \displaystyle  q 
\sum_{j=0}^{u-1} 
(\eta + q)^{*(u-1-j)}\theta^j/j! & & 
\mbox{for $d=g$;} \\
\rm{(iii)} & &
(\eta + r)^{*u} & - & \displaystyle  r 
\sum_{j=0}^{u-1} 
(\eta + r)^{*(u-1-j)}\theta^j/(j+1)! & &
\mbox{for $d=g-1$, 
where $r=q \, \theta/(1-q)$;} \\
\rm{(iv)} & &
\eta^{*u} & - & \displaystyle  q \sum_{j=0}^{u-1} \eta^{*j}
\theta^{n-j}/(n-j)! & \\ 
& & & 
\setlength\arraycolsep{2pt}
\renewcommand{\arraystretch}{1.2}
\begin{array}{c} + \\ \phantom{+} \end{array}
& 
\setlength\arraycolsep{2pt}
\renewcommand{\arraystretch}{1.2}
\begin{array}{l}
q \, u \, \eta^{*(u-1)}\theta^{g-d}/(g-d)! + O(q^2) \\
\phantom{xxx} \end{array}
& & 
\setlength\arraycolsep{2pt}
\renewcommand{\arraystretch}{1.2}
\begin{array}{l}
\mbox{for $g/2 < d \leq g-1$,} \\
\phantom{xxx}\mbox{where $n=u-d+g-1$.} \end{array}  
\end{array}$$
\es

\pf.  For $d>g$, by \re{aa}(i)
$$\eta * \eta^v = \eta^{v+1} + q \, \frac{\theta^{g-d+v}}{(g-d+v)!}.$$
Part (i) follows inductively by taking the quantum product of both
sides with $\eta$.  Similarly, for $d=g$
$$(\eta + q) * \eta^v = \eta^{v+1} + q \, \theta^v/v!,$$ 
and for $d=g-1$
$$(\eta + r) * \eta^v = \eta^{v+1} + r \, \theta^v/(v+1)!,$$ 
so parts (ii) and (iii) also follow inductively.  Part (iv) follows
from \re{t} by an easy induction.  \fp


For $d \geq g-1$ or $d < \ratio{4}{5}g-\ratio{3}{5}$, it is now easy
to obtain a complete set of quantum relations which are explicit, if
somewhat inelegant.  Just plug the formulas of \re{w} into \re{v}.
When $\ratio{2}{3}g+\ratio{1}{3} \leq d < \ratio{4}{5}g -
\ratio{3}{5}$, terms of higher order in $q$ appear in the quantum
product.  However, these do not affect the relations, since by the
last statement in \re{b}, the latter are homogeneous of degree $d$ or
$d+1$, too large to contain a power of $q^2$.  The same is true for
the linear terms in $q$ when $\half g + 1 \leq d < \ratio{2}{3}
g-\ratio{1}{3}$, leading to the amusing phenomenon that in that range,
the quantum and classical rings are isomorphic, but only by a map
which is not the identity.

Although the general quantum relation is no thing of beauty, there are
a few exceptions.

{\samepage \bs{Corollary}
\label{y}
\bl 
\item For $d > 2g-2$,
$$ \eta^{*(d-2g+1)} *
\prod_{i=1}^g \! \raisebox{1.5ex}{$*$} 
(\eta - \si_i) = q;$$
\item For $g < d \leq 2g-2$,
$$\prod_{i=1}^g \! \raisebox{1.5ex}{$*$} 
(\eta - \si_i) = q \, \eta^{*(2g-1-d)};$$
\item For $d=g$,
$$\prod_{i=1}^g \! \raisebox{1.5ex}{$*$} 
(\eta - \si_i + q) = q \, (\eta+q)^{*(g-1)}.$$
\el
\es}

\pf.  The first relation follows directly from \re{c} and \re{w}(i),
since the coefficient of $q$ in the latter is $-1$ for $u=d-g+1$ and $0$
for $u<d-g+1$.  To prove (ii) and (iii), note first that by \re{b}, the
relation
$$0 = \prod_{i=1}^g (\eta-\si_i) $$
holds in the classical ring for $d \leq 2g-2$.  This is equivalent to 
$$0 = \sum_{j=0}^g (-1)^j \eta^{g-j} \theta^j/j!.$$
One plugs in the formulas of \re{w}(i) and (ii) for $\eta^u$ and notes
that the double sum in the second term telescopes using the binomial
theorem. \fp

\bit{Relation with Givental's work}
In fact, relations (ii) and (iii) in the corollary above were known to
the authors before any of the Gromov-Witten invariants, and were the
starting point of the investigation.  With a little ingenuity, they
can be read off from formulas in the wonderful paper of Givental
\cite{g}, which is concerned chiefly with proving that the
Gromov-Witten invariants of the quintic threefold satisfy the
Picard-Fuchs equation.  We will sketch an outline of the connection.

At the heart of Givental's paper is a ``quantum Lefschetz hyperplane
theorem'' explaining how the Gromov-Witten invariants of a variety are
related to those of a hyperplane section.  This allows him to compute
a generating function for the Gromov-Witten invariants of complete
intersections.  He shows in Corollary 6.4 that any differential
operator annihilating this generating function determines a relation
in the quantum cohomology.

The symmetric products of a curve do not appear in any natural way as
complete intersections in a projective space.  But they do appear as
complete intersections in a projective {\em bundle}.  Indeed, for any
$d \leq 2g-2$, choose a reduced divisor $P = \sum p_i$ of degree
$2g-1-d$.  Then $D \mapsto D+P$ gives an embedding $\Si_d
\hookrightarrow \Si_{2g-1}$, whose image is a complete intersection of
divisors in the linear system of $\co(1)$.  In particular, for $d \leq
2g-2$, $\Si_d$ embeds in $\Si_{2g-1}$, which is a $\Pj^{g-1}$-bundle
over $\Jac_{2g-1}$.

Givental's methods are perfectly adapted to studying this more general
case.  Indeed, his formulas for Gromov-Witten invariants are derived
as special cases of formulas for equivariant Gromov-Witten
invariants.  These can be regarded as universal formulas for
Gromov-Witten invariants of complete intersections in projective
bundles, because they are defined as relative Gromov-Witten invariants
for complete intersections in the universal projective bundle over
the classifying space.  

Givental works with the group $(\C^\times)^{n+1}$ acting on $\Pj^n$.
The classifying space is then $(\C\Pj^{\infty})^{n+1}$, and the
universal bundle is a direct sum of line bundles.  This would appear
to be a problem, because $\Si_{2g-1}$ is not the projectivization of a
direct sum of line bundles.  However, a splitting principle argument
shows that all of Givental's equivariant formulas extend word for word
to the action of $\GL{n,\C}$.  The classifying space is then an
infinite Grassmannian; and any projective bundle at all is pulled back
from some map to a Grassmannian, even in the algebraic category.

Thus a formula in equivariant cohomology determines a formula in the
cohomology of $\Si_{2g-1}$, or more properly, in the cohomology of the
spaces of stable maps to $\Si_{2g-1}$.  This is perhaps surprising, since
no group is acting on $\Si_{2g-1}$.  However, the group is present
as the maximal torus of the structure group of the projective bundle
$\Si_{2g-1} \to \Jac_{2g-1}$; one regards equivariant cohomology as
giving universal formulas in the cohomology of all such bundles.

Another apparent problem is that the equivariant methods treat only
Gromov-Witten invariants for classes $e \in H_2$ killed by the
projection to the base, namely $\Jac_{2g-1}$.  But luckily, as we have
seen, these are the only nonzero invariants of $\Si_d$.

The symmetric product $\Si_d$ is an equivariant complete intersection
in $\Si_{2g-1}$ of type $(1,1,\dots, 1)$.  Such an intersection in a
single projective space would just be a linear subspace, but it can be
less trivial equivariantly, that is to say, in families.  Givental's
results therefore apply with $r=2g-1-d$, $l_1 = \cdots = l_r = 1$.

The case $d>g$ is covered by Givental's \S9.  The equivariant
quantum potential satisfies the differential equation shown below his
Theorem 9.5.  As instructed in Corollary 6.4, we make the
substitutions $\hbar \, d/dt = \eta$ (the relative hyperplane class),
$e^t = q$, $\hbar = 0$ to obtain a quantum relation.  We also
substitute $\la_i = \si_i$ since by \re{c} these are the Chern roots of
the bundle whose projectivization is $\Si_{2g-1}$, and $\la'_i = 0$
since these come from an additional equivariance we are not using.
The result is precisely \re{y}(ii)!

Similarly, the case $d=g$ is covered by \S10.  The differential
equation satisfied by the equivariant quantum potential $S'$ is not
explicitly stated, but it is
$$\prod_i (D-\la_i) \, S' 
= l_1 \cdots l_r e^t \prod_{j,m} (l_j D - \la'_j + m \hbar) \, S',$$
where $D = \hbar \, d/dt + l_1 \cdots l_r e^t$ as in Corollary 10.8.
Making the same substitutions as before transforms $D$ to $\eta + q$
and this equation to \re{y}(iii).

Finally, the case $d=g-1$ is analogous to the Calabi-Yau case, which
is covered by \S11.  Again the differential operator determines a
relation in degree $2g$, but now this is vacuous, since the real
dimension of $\Si_{g-1}$ is only $2g-2$.  This is a familiar
occurrence in Givental's work: the Picard-Fuchs equation, for
example, has degree 4, so it gives a quantum relation in $H^8$ of the
quintic threefold, which is of course trivial.  

Nevertheless, the Picard-Fuchs equation on the quintic does contain
valuable information: enough to determine the 3-point invariants at
genus 0, and hence the virtual number of rational curves of all
degrees.  On $\Si_{g-1}$ the corresponding invariants were all worked
out in \re{ii}, and are all determined by the lines using the
Aspinwall-Morrison formula.  It is like a Calabi-Yau with no higher
degree rational curves.  But it certainly ought to be possible to
recover these invariants for $\Si_{g-1}$, or even $\Si_d$ for $d \geq
g-1$, by calculating Givental's quantum potential in this case.

A more daunting project would be to extend these methods to $\Si_d$
for $d < g-1$, where the results of this paper are incomplete.
Givental remarks that the corresponding case for complete
intersections remains unsolved, but is in a sense ``less
interesting,'' since nonzero invariants appear only in finitely many
degrees.  For symmetric products, however, this is clearly the most
interesting and mysterious case.

\bit{The homotopy groups \\ of symmetric products of curves} 
In this appendix, we first compute the fundamental group of a
symmetric product, then show that the Hurewicz homomorphism $\pi_2 \to
H_2$ has rank 1 for $d>1$.  These results are not needed anywhere
else, but they clarify the non-contribution of other homology classes,
from the point of view of symplectic topology.

For $d>2g-2$, the symmetric product is a projective bundle over the
Jacobian, so it has fundamental group $H_1(\Si;\Z)$.  On the other hand,
for $d=1$ it is of course $\pi_1(\Si)$. What happens in between?  The
simplest possible thing, it turns out.  We are grateful to Michael
Roth for supplying the following theorem and its proof; but see also
Grothendieck \cite{gr}.

\bs{Theorem}  For $d>1$, $\pi_1(\Si_d) = H_1(\Si; \Z)$.  Indeed,
the Abel-Jacobi map induces an isomorphism on fundamental groups.
\es

\pf.  Choose a basepoint $p \in \Si$ and let $i: \Si \to \Si^d$ be
given by $x \mapsto (x,p,\dots, p)$.  Also let $\psi: \Si^d \to \Si_d$ be
the quotient by the symmetric group.  We will show first, that
$\pi_1(\psi i)$ surjects on $\pi_1(\Si_d)$, second, that its image is
abelian, and third, that the kernel is exactly the commutator
subgroup.

Choose as a basepoint in $\Si_d$ the divisor $d \cdot p$.  Since the
big diagonal has real codimension 2, any loop can be perturbed so that
it meets the big diagonal only at the basepoint.  It then lifts
unambiguously to a loop on $\Si^d$.  Hence $\pi_1(\psi)$ is 
surjective. But any loop on $\Si^d$ is homotopic to a composition of
loops on the various factors.  Since $\psi$ is symmetric, each of these
loops can be replaced with the corresponding loop on the first factor,
that is, on the image of $i$, without changing $\pi_1(\psi)$ of the
composite.  Hence $\pi_1(\psi i)$ is surjective.

Let $\ga_2, \ga_2$ be two loops in $\pi_1(\Si)$.  To show that their
images in $\Si_d$ commute, first include them in $\Si^d$ via $i$.
Then note that $\ga_2$ can be transferred to the second factor without
changing its image in $\Si_d$; then, however, it commutes with $\ga_1$
in $\pi_1(\Si^d) = \pi_1(\Si)^d$.

Finally, to show that $\psi i$ kills only loops in the commutator,
compose $\psi i$ with the Abel-Jacobi map.  The resulting map $\Si \to
\Jac_d$ is itself just the composition of the Abel-Jacobi map with the
identification $\Jac_1 \cong \Jac_d$ induced by $p$.  But it is
well-known that the Abel-Jacobi map induces an isomorphism $H_1(\Si;
\Z) = H_1(\Jac_1; \Z)$; see for example Griffiths-Harris
\cite[2.7]{gh}. \fp

\bs{Theorem}
For $d>1$, the Hurewicz homomorphism $\pi_2(\Si_d) \to H_2(\Si; \Z)$
has rank $1$.
\es

\pf. It is easy to see that the rank is at least $1$: for example,
take $\Si$ to be hyperelliptic; there is then a one-dimensional linear
system $\ell$ on $\Si_d$, which is a 2-sphere, non-trivial in rational
homology because $\eta \cdot \ell = 1$.

To show it is at most $1$, note that by the above theorem the
universal cover $\tilde{\Si}_d$ is the fibered product $\Si_d
\times_{\Jac_d} H^1(\Si, \co)$, where the vector space $H^1(\Si, \co)$
is the universal cover of $\Jac_d$.  Now the push-forward $H_2(\Si_d;
\Q) \to H_2(\Jac_d; \Q)$ has $1$-dimensional kernel: this follows for
example from \re{x}(iii) and \re{b}.  Since the map $\tilde{\Si}_d \to
\Si_d \to \Jac_d$ also factors through the contractible space
$H^1(\Si, \co)$, the induced map on $H_2$ is zero and so the map
$H_2(\tilde{\Si}_d) \to H_2(\Si_d)$ has rank at most $1$.  But by the
Hurewicz isomorphism and the homotopy exact sequence,
$H_2(\tilde{\Si}_d) = \pi_2(\tilde{\Si}_d) = \pi_2(\Si_d)$, and the
map to $H_2(\Si_d)$ is exactly the Hurewicz homomorphism.  \fp

\end{document}